\documentclass[12pt]{article}

\usepackage{multirow}
\usepackage{mathtools}
\usepackage{subcaption}
\usepackage{amsfonts,amssymb,amsmath,amsthm,latexsym,bbm}
\usepackage{enumerate}
\usepackage{epsf,epsfig}
\usepackage{xcolor,colortbl,color}
\usepackage{graphicx,graphics}
\usepackage{caption}
\usepackage[utf8]{inputenc}
\usepackage{tikz}
\usepackage{stmaryrd }
\usepackage{epstopdf}

\usepackage{comment}
\usepackage{url}
\makeatletter \@addtoreset{equation}{section} \makeatother
\makeatletter \@addtoreset{enunciato}{section} \makeatother

\renewcommand{\theequation}{\thesection.\arabic{equation}}
\newcommand{\be}[1]{\begin{equation}\label{#1}}
\newcommand{\ee}{\end{equation}}
\newcommand{\bl}[1]{\begin{lemma}\label{#1}}
\newcommand{\el}{\end{lemma}}
\newcommand{\br}[1]{\begin{remark}\label{#1}}
\newcommand{\er}{\end{remark}}
\newcommand{\bt}[1]{\begin{theorem}\label{#1}}
\newcommand{\et}{\end{theorem}}
\newcommand{\bd}[1]{\begin{definition}\label{#1}}
\newcommand{\ed}{\end{definition}}
\newcommand{\bp}[1]{\begin{proposition}\label{#1}}
\newcommand{\ep}{\end{proposition}}
\newcommand{\bc}[1]{\begin{corollary}\label{#1}}

\newcommand{\bcj}[1]{\begin{conjecture}\label{#1}}
\newcommand{\ecj}{\end{conjecture}}

\def\<{\langle}
\def\>{\rangle}

\def\eatspace#1{\relax}
\def\unskipit{\expandafter\eatspace}

\def\claim#1{\begin{trivlist}\item[\hskip\labelsep\bf#1]\it}
\def\endclaim{\end{trivlist}}

\numberwithin{equation}{section}

\newtheorem{theorem}{Theorem}[section]
\newtheorem{lemma}[theorem]{Lemma}
\newtheorem{pr}{Proposition}[section]
\newtheorem{corollary}[theorem]{Corollary}
\newtheorem{definition}{Definition}[section]
\newtheorem{remark}[theorem]{Remark}

\newcommand{\eproof}{{\mbox{\ }~\hfill
\mbox{\large $\Box$} \par \vskip 10pt}}
\newcommand{\pf}{\noindent{\bf Proof}}

\title{Classical unique continuation property for multi-terms time fractional diffusion equations}
\author{Ching-Lung Lin\thanks{Department of Mathematics, National Cheng-
Kung University, Tainan 701, Taiwan.  (Email:
cllin2@mail.ncku.edu.tw)}\qquad Gen Nakamura\thanks{Department of
Mathematics, Hokkaido University, Sapporo 060-0808, Japan and Research Center of Mathematics for Social Creativity, Research Institute for Electronic Science, Hokkaido Uni-
versity, Sapporo, 060-0811, Japan.
\newline
(Email: gnaka@math.sci.hokudai.ac.jp)
}}

\date{}

\begin{document}
\renewcommand{\theequation}{\thesection.\arabic{equation}}
 \maketitle
\begin{abstract}
As for the unique continuation property (UCP) of solutions in $(0,T)\times\Omega$
with a domain $\Omega\subset{\mathbb R}^n,\,n\in{\mathbb N}$ for a multi-terms time fractional diffusion equation,
we have already shown it by assuming that the solutions are zero for $t\le0$ (see \cite{LN2019}).
Here the strongly elliptic operator for this diffusion equation
can depend on time and the orders of its time fractional derivatives are in $(0,2)$.
This paper is a continuation of the previous study. The aim of this paper is
to drop the assumption that the solutions are zero for $t\le0$.
We have achieved this aim by first using the usual Holmgren transformation
together with the argument in \cite{LN2019} to derive the UCP in $(T_1,T_2)\times B_1$
for some $0<T_1<T_2<T$ and a ball $B_1\subset\Omega$. Then if $u$ is the solution of the equation
with $u=0$ in $(T_1,T_2)\times B_1$, we show $u=0$ also in
$((0,T_1]\cup[T_2,T))\times B_r$ for some $r<1$ by using the argument in \cite{LN2019}
which uses two Holmgren type transformations different from the usual one.
This together with spatial coordinates transformation, we can obtain the usual UCP
which we call it the classical UCP given in the title of this paper for
our time fractional diffusion equation.
\end{abstract}

\section{Introduction}\label{sec1}
We first define an anomalous diffusion equation with multi-terms fractional time derivatives. As mentioned in Abstract, this paper is a continuation of our previous study on the unique continuation property abbreviated by UCP of solutions of the anomalous diffusion equation given in the paper \cite{LN2019}, we will basically use the same notations as those in that paper for the readers convenience. To begin with, let $(0,T)$ be a time interval and $\Omega\subset{\mathbb R}^n$ with $n\in{\mathbb N}$ be a domain.
 Also we denote by $\emph{L}$ the strongly elliptic operator defined by
\begin{equation}\label{L}
\emph{L} u(t,y):=\sum_{j,k=1}^na_{jk}(t,y)\partial_{y_j}\partial_{y_k}u(t,y)
\end{equation}
with $a_{jk}(t,y)=a_{kj}(t,y)\in C^\infty([0,T]\times\overline\Omega),\,1\le j,k\le n$ such that
\begin{equation}\label{1.2}
\begin{array}{l}
\delta |\xi|^2 \leq \sum_{j,k=1}^na_{jk}(t,y)\xi_j\xi_k \leq \delta^{-1} |\xi|^2
\end{array}
\end{equation}
holds for $(t,y)\in [0,T]\times\overline\Omega,\, \xi=(\xi_1,\cdots,\xi_n)\in{\mathbb R}^n$ with some constant $\delta>0$. Henceforth we always extend each $a_{jk}\in C^\infty([0,T]\times\overline\Omega)$ to
$a_{jk}\in B^\infty({\mathbb R}^{1+n})$ without destroying \eqref{1.2} and $a_{jk}(t,y)=a_{kj}(t,y)$, where
$B^\infty({\mathbb R}^{1+n})$ denotes the set of all functions in $C^\infty({\mathbb R}^{1+n})$ which
are bounded together with their derivatives. To avoid any confusion, we note that the notation $A\subset B$ describing $A$ is a subset of $B$ admits the case $A=B$.

Now let $m\in {\mathbb N}$ and $\alpha_j, q_j,$ $j=1,2\cdots m$ be constants such that $q_1=1$ and $2>\alpha=\alpha_1>\alpha_2>\cdots>\alpha_m>0$ and consider an anomalous diffusion equation with multi-terms fractional time derivatives given by
\begin{equation}\label{anomalous eq}
\Sigma_{j=1}^mq_j\partial_t^{\alpha_j} u(t,y)-\emph{L} u(t,y)=\ell_1(t,y;\nabla_y)u(t,y),\,\,t\in(0,T),\,y\in\Omega,
\end{equation}
where
$\ell_1(t,y;\nabla_y)$ is a linear partial differential operator of order one with $C^\infty([0,T]\times\overline\Omega)$ coefficients and $\partial_t^\alpha u$ is defined as follows.
\begin{equation}\label{1.4}
\partial_t^\alpha u(t,y)=
\left\{
\begin{array}{ll}
\frac{1}{\Gamma(k-\alpha)}\big((t^{k-1-\alpha}H(t)\otimes\delta_y)*\partial^k_t u\big)(t,y)\,\,&\text{for}\,\,k-1<\alpha<k,\,k=1,2,\\\\
\partial_t u(t,y)\,\,&\text{for}\,\,\alpha=1
\end{array}
\right.
\end{equation}
with the Heaviside function $H(t)$ and the Dirac delta function $\delta_y$ supported at $y$.
Hence $\partial^\alpha$ for $0<\alpha<1$ or $1<\alpha<2$ is the Caputo fractional derivative of order $\alpha$. The coefficients of $\ell_1(t,y;\nabla_y)$ could be just in $C^0([0,T]\times\overline\Omega)$, but to simplify our arguments we have assumed that they are in $C^\infty([0,T]\times\overline\Omega)$.
We remark here whenever $u$ satisfies \eqref{anomalous eq},
we are assuming that the Caputo fractional derivatives $\partial_t^{\alpha_j}u,\,1\le j\le m$ can be defined.

Before giving our main result, we first would like to clarify the novelty of our main result by giving our previous result of the UCP which was proven by using the Holmgren type transformation different from the usual Holmgren transformation.

\begin{theorem}\label{thm1.0}{\rm (our previous UCP)}
Let $u\in H^{\alpha,2}((0,T)\times\Omega)$
(see below in this section for its definition) be a solution of \eqref{anomalous eq} supported on $\{t\ge0\}$.
Then, $u=0$ in a subdomain of $\Omega$ over $(0,T)$ implies
$u=0$ in $(0,T)\times\Omega$. Here the assumption that $u$ is supported on $\{t\ge 0\}$ means that its zero extension across $t=0$ belongs to $H^{\alpha,2}((-\infty,T)\times\Omega)$.
We quote this as {\rm the initial condition}.
\end{theorem}

Now we are ready to give our main result as follows.

\begin{theorem}\label{thm1.1}{\rm(classical UCP)}
Let $u\in H^{\alpha,2}((0,T)\times\Omega)$
be a solution of \eqref{anomalous eq}.
Then, $u=0$ in a subdomain of $\Omega$ over $(0,T)$ implies
$u=0$ in $(0,T)\times\Omega$.
\end{theorem}

\begin{remark}
This is the usual UCP which we call it the classical UCP given in the title of this paper for solutions of \eqref{anomalous eq}. To be consistent with Theorem \ref{thm1.0}, we have assumed that $q_l,\,1\le l\le m$ are constants. But Theorem \ref{thm1.0} and Theorem \ref{thm1.1} still hold even for the case $q_l=q_l(x)\in C^\infty(\overline\Omega),$ $1\le l\le m$ with $q_1\equiv 1$, and their proofs remain the same as those of Theorem \ref{thm1.0} and Theorem \ref{thm1.1}, respectively.
\end{remark}

Theorem \ref{thm1.1} follows from the following theorem together with spatial coordinates transformation.

\begin{theorem}\label{thm1.3}
Let $B_1\subset\Omega$ and $u\in H^{\alpha,2}((0,T)\times\Omega)$ be a solution of \eqref{anomalous eq}.
If $u=0$ in  $(0,T)\times B_r$ for some $0<r<1$ and $z\in \partial B_r$, then there exists a small $0<\hat{r}<r$ such that
$u=0$ in $(0,T)\times B_{\hat{r}}(z)$.
Here $B_r, B_r(z)$ are defined as
\begin{equation*}
B_r=\{y=(y_1,\cdots,y_n): |y|<r\},\,\,B_r(z)=\{y=(y_1,\cdots,y_n):y-z\in B_r \}
\end{equation*}
with $|y|=\sqrt{\sum_{j=1}^n y_j^2}$ for $y=(y_1,\cdots,y_n)$.
\end{theorem}

\begin{remark}
This result will be proven by using both the usual Holmgren transformation, the aforementioned Holmgren type transformation and another Holmgren type transformation (see \eqref{H type trf}).
For the simplicity of notations, we just took $B_1\subset\Omega$ as a ball centered at the origin with radius $1$.
\end{remark}

\medskip
 The anomalous diffusion equation was used to model for examples, the anomalous diffusion phenomena observed in a highly heterogeneous acquifer (\cite{Adams}, \cite{hat1} and the reference there in) and a complex viscoelastic material (\cite{Brown}, \cite{Ginoa} and the references therein). Then, in the past decades, it has become more and more popular in interdisciplinary fields due to its flexibility of modeling various nonlocal phenomena.
 \par
 The classical UCP for the anomalous diffusion equation was obtained for some special cases by some Carleman estimates. More precisely, it was given for the one space dimensional case with one time fractional derivative of order $1/2$ in \cite{XCY}, and its generalization to the general space dimensional case was given in \cite{CLN}. For the case that the order of the time fractional derivative is in $(0,1)$, the UCP under the aforementioned initial condition was given in \cite{LN} and its generalization to the muliti-terms time fractional derivative case with orders in $(0,2)$ was given in \cite{LN2019}.

Besides the classical UCP, there is another UCP called the weak UCP which is for any solution $u$ of an initial boundary value problem for anomalous diffusion equation with initial data at $t=0$ and homogeneous boundary condition. This weak UCP states that if this solution satisfies $u=0$ in an open subset of $(0,T)\times\Omega$, then $u=0$ in $(0,T)\times\Omega$. For the fractional anomalous diffusion equation with single term time fractional derivative, it was proved in \cite{Li} for the one space dimensional case and it was also proved in \cite{Jiang} for the higher space dimensional case. Even a further generalization of \cite{Jiang} has been done in \cite{LHY} for the case that the order of the time fractional derivative is in $(0,2)$ and the elliptic part of the equation is the minus Laplacian. Furthermore, the weak UCP was used to solve the uniqueness of some inverse source problems (\cite{Jiang}, \cite{LHY}, \cite{Sun}  and the references there in). This already revealed the importance of the classical UCP in the studies of many inverse problems for the anomalous diffusion equations. Not only for the inverse problems but also for the control problems, one of the very important key in studying these problems is the classical UCP, and it has been awaited very much for a long time.
As usual we will derive some Carleman estimates for solutions of the anomalous diffusion equation to prove the usual UCP. More precisely, as mentioned before in Abstract, we can have a Carleman estimate associated to the usual Holmgren transformation and another two Carleman estimates associated to the two Holmgren type transformations to prove Theorem \ref{thm1.3}. The basic idea behind the derivation of Carleman estimates is as follows. Likewise the usual diffusion equation, the anomalous diffusion equation is semi-elliptic. Based on this observation, we will use Treve's argument (\cite{Treves} and its nice introduction given in \cite{Taylor}) for the derivation of Carleman estimates.

The rest of this paper is devoted to proving Theorem \ref{thm1.3} and it is organized as follows.
In Section 2, we state Proposition \ref{prop2.1} and
give its proof. By this proposition, we can have the UCP in $(T_1,T_2)\times B_{\tilde r}(z)$ for some $0<T_1<T_2<T$
and a ball $B_{\tilde r}(z)\subset\Omega$ by using the usual Holmgren transformation.
A very brief outline of its proof is as follows.
We first observe that the principal symbol of the anomalous diffusion operator conjugated by $e^t$, also undergone the usual Holmgren transformation and replace $D_x$ by $D_x+i|D_z|\nabla\psi$ with the Carleman weight function $\psi=(x_n-2X)^2/2$ by adding a new variable $z$ (see \eqref{eq:3.3}) coincides with that of obtained by using the Holmgren type transformation.
Hence just using the same argument given in \cite{LN2019},
we can derive a Carleman estimate which yields Proposition \ref{prop2.1}.
In Section 3, we will prove Theorem \ref{thm1.3}. Its very brief outline of the proof is as follows.
In Proposition \ref{Pro3.1}, we will show that we can have $u=0$ in $[T_1,T)\times B_{\tilde r}(z)$ by using Theorem \ref{thm1.0} and Proposition \ref{prop2.1}. Then,
in Proposition \ref{Pro3.2}, we will show $u=0$ in $(0,T_1]\times B_{\hat{r}}(z)$ with $\hat{r}<\tilde{r}$ in Section 3.
The argument for showing this is done in the same way as given in \cite{LN2019}
except using the following another Holmgren type transformation \begin{equation}\label{H type trf}
x'=\tilde y',\,x_{s,n}=\tilde y_n+c|\tilde y'|^2-(\frac{sX}{T}t-X),\,\tilde t=t,
\end{equation}
where $\tilde y=(\tilde y',\tilde y_n)=(\tilde y_1,\cdots,\tilde y_{n-1},\tilde y_n)$, $|\tilde y'|=\sqrt{\sum_{j=1}^{n-1}\tilde y_j^2}$.

Finally we close this section by giving the definition of Sobolev space
$H^{\alpha,2}({{\mathbb R}}\times\Omega)$. For $m,s\in {\mathbb R}$,
$v=v(t,x)\in \mathcal{S}'({\mathbb R}_t^1\times{\mathbb R}_x^n)$ belongs to the function space
$H^{m,s}({\mathbb R}_t\times{\mathbb R}_x^n)$ if
$$\|v\|^2_{H^{m,s}}:=\iint(1+|\xi|^s+|\tau|^m)^2|\hat{v}|^2d\tau d\xi$$
is finite, and $\|v\|^2_{H^{m,s}}$ denotes the norm of $v\in H^{m,s}({\mathbb R}_t\times{\mathbb R}_x^n)$, where
$\hat{v}$ is the Fourier transform defined by
$$\mathcal{F}(v)(\tau,\xi)=\hat{v}(\tau,\xi)=\iint e^{-it\tau-ix\cdot\xi}v(t,x)dxdt$$
and $\mathcal{S}'({\mathbb R}_t^1\times{\mathbb R}_x^n)$ is the dual space
of the space $\mathcal{S}({\mathbb R}_t^1\times{\mathbb R}_x^n)$ of rapidly decreasing functions in ${\mathbb R}_t^1\times{\mathbb R}_x^n$.
Further, for any open sets $A\subset{{\mathbb R}}_t,\,B\subset{{\mathbb R}}_x^n$,
we define $H^{m,s}(A\times B)$ as the restriction of $H^{m,s}({\mathbb R}_t\times{\mathbb R}_x^n)$ to $A\times B$.

\section{A local uniqueness result}\label{sec2}
\setcounter{equation}{0}
This section is devoted to proving the following local uniqueness result.

\begin{pr}\label{prop2.1}
Let $u\in H^{\alpha,2}((0,T)\times\Omega)$ be a solution of \eqref{anomalous eq}.
Suppose $u=0$ in  $(0,T)\times B_r$ for some $0<r<1$. Take any $z\in\partial B_r$.
Then, there exists a small $0<\tilde{r}<r$ and $0<T_1<T_2<T$ such that
$u=0$ in $(T_1,T_2)\times B_{\tilde{r}}(z)$.
\end{pr}
\pf${}\,$
The proof of Proposition \ref{prop2.1} is quite similar to that of Theorem 1.2 in \cite{LN2019}.
The only difference here is that we use the usual Holmgren transformation.
It is important to note here that the alternatives to the Holmgren type transformation and the initial condition
which we used to prove Theorem 1.2 in \cite{LN2019} is the usual Holmgren transformation
and some cutoff functions attached to this transformation. To avoid any misunderstanding,
let us emphasize that we do not need any initial condition for the proof.
We will only give the steps to prove the theorem
clarifying the role of the usual Holmgren transformation and the cutoff functions.

By a translation and a rotation of coordinates, let $z$ be given in the form $z=(z_1,\cdots,z_{n-1},0)$. To prove Proposition \ref{prop2.1},
we use the usual Holmgren transformation
\begin{equation}\label{2.1}
\begin{array}{l}
x'=y'-z', x_n=y_n+c|y'-z'|^2+c(t-\frac{T}{2})^2, \tilde{t}=t
\end{array}
\end{equation}
with $y'=(y_1,\cdots,y_{n-1})$, $z'=(z_1,\cdots,z_{n-1})$ and $x'=(x_1,\cdots,x_{n-1})$, where $c\geq 1$ is a fixed positive constant.
We will quote this new coordinates $(t,x)=(\tilde t,x)$ as the {\sl Holmgren coordinates}.

By $\partial_t=\partial_{\tilde{t}}+2c(\tilde{t}-\frac{T}{2})\partial_{x_n}$,
$ \partial_{y_j}=\partial_{x_j}+2cx_j\partial_{x_n}$ for $j=1,\cdots,n-1$ and $\partial_{y_n}=\partial_{x_n}$, we have
\begin{equation}\label{2.2}
\begin{array}{l}
\partial_{t}^\alpha u(t,y)=\frac{1}{\Gamma(1-\alpha)}\int_0^t(t-\eta)^{-\alpha}\partial_\eta u(\eta,y)d\eta\\
=\frac{1}{\Gamma(1-\alpha)}\int_0^{\tilde{t}}(\tilde{t}-\tilde{\eta})^{-\alpha}\partial_{\tilde{\eta}} u(\tilde{\eta},x)d\tilde{\eta}
+\frac{2c(\tilde{t}-\frac{T}{2})}{\Gamma(1-\alpha)}\int_0^{\tilde{t}}(\tilde{t}-\tilde{\eta})^{-\alpha}\partial_{x_n}u(\tilde{\eta},x)d\tilde{\eta}
\end{array}
\end{equation}
for $0<\alpha<1$ and
\begin{equation}\label{another 2.2}
\begin{array}{ll}
\partial_{t}^\alpha u(t,y)=\frac{1}{\Gamma(2-\alpha)}\int_0^t(t-\eta)^{1-\alpha}\partial_\eta^2 u(\eta,y)d\eta\\
=\frac{1}{\Gamma(2-\alpha)}\int_0^{\tilde{t}}(\tilde{t}-\tilde{\eta})^{1-\alpha}\partial_{\tilde{\eta}}^2 u(\tilde{\eta},x)d\tilde{\eta}
+\frac{4c(\tilde{t}-\frac{T}{2})}{\Gamma(2-\alpha)}\int_0^{\tilde{t}}(\tilde{t}-\tilde{\eta})^{1-\alpha}\partial_{\tilde\eta}\partial_{x_n}u(\tilde{\eta},x)d\tilde{\eta}\\
\quad+\frac{1}{\Gamma(2-\alpha)}\int_0^{\tilde t}(\tilde t-\tilde\eta)^{1-\alpha}
\big(4c^2(\tilde\eta-\frac{T}{2})^2\partial_{x_n}^2+2c\partial_{x_n}\big)u(\tilde\eta,x)\,d\tilde\eta
\end{array}
\end{equation}
for $1<\alpha<2$, where $\Gamma(\cdot)$ denotes the Gamma function and $u(\tilde{t},x)$
with $\tilde{t}=\tilde{\eta}$ on the right hand side of \eqref{2.2}, \eqref{another 2.2}
is the push forward of $u(t,y)$ by the change of variables.

It should be noted that
$ x_n=y_n+c|y'-z'|^2+c(t-\frac{T}{2})^2$, and $u(\tilde{t},x)$ satisfies
$${\rm supp}\,  u\subset\{x_n\geq 0\}.$$
Further, let $\chi$ be a smooth function such that
\begin{equation}\label{2.3}
\chi (x_n)=
\begin{cases}
\begin{array}{l}
1,\quad x_n\leq (1-\epsilon)X,\\
0,\quad x_n\geq X.
\end{array}
\end{cases}
\end{equation}
Then
\begin{equation}\label{2.4}
{\rm supp}\, \chi u\subset \{0\leq x_n\leq X\}.
\end{equation}
Also, we will assume that
\begin{equation}\label{2.5}
    |x'|\lesssim \sqrt{X}, \quad x_n\leq X,
\end{equation}
here and after $A \lesssim B$ denotes $A \leq CB$ with positive constant $C$ depending on $n$ and $\alpha$.

To avoid having the branch points for the symbols of $\partial_t^{\alpha_l},\,1\le l\le m$, we conjugate $\sum_{l=1}^m\,q_l\,\partial_t^{\alpha_l}-\emph{L}$ with
$e^{t}$. Thus we consider the operator $e^{-t}(\sum_{l=1}^m\,q_l\,\partial_t^{\alpha_l}-\emph{L})e^{t}$ and denote it by $P(t,x,D_t,D_x)$
in terms of the Holmgren coordinates $(t,x)$ (see \eqref{2.1}), where $D_t=-i\partial_t$, $D_x=-i\partial_x$ with $i=\sqrt{-1}$. Then it is easy to see that for example for the case $0<\alpha_\ell<1,\,1\le\ell\le m$, the total symbol
$p(t,x,\tau,\xi)$ of $P(t,x,D_t,D_x)$ is given as
\begin{equation}\label{2.6}
\begin{array}{rl}
p(t,x,\tau,\xi)=&\Sigma_{l=1}^m q_l\,(1+i\tau)^{\alpha_l}+a_{nn}(t,x)\xi_n^2+2\Sigma_{j=1}^{n-1}a_{jn}(t,x)\xi_n(\xi_j+2cx_j\xi_n)\\
&+\Sigma_{j,k=1}^{n-1}a_{jk}(t,x)(\xi_j+2cx_j\xi_n)(\xi_k+2cx_k\xi_n)\\
&+\Sigma_{l=1}^m\,q_l\,2c(t-\frac{T}{2})i^\alpha(\tau-i)^{\alpha_l-1}\xi_n.
\end{array}
\end{equation}

Now let $\psi=\frac{1}{2}(x_n-2X)^2$ be the Carleman weight function and consider the symbol $p(t,x,\tau,\xi+i|\sigma|\nabla\psi)$ over $\mathbb{R}^{n+1}\times\mathbb{R}_z$ to define a pseudo-differential operator $P_\psi(t,x,D_t,D_{x},D_z)$ by
$$P_\psi=P_\psi(t,x,D_t,D_{x},D_z)=p(t,x,D_{t},D_x+i|D_z|\nabla \psi)$$
which is given for any compactly supported distribution $v$ in $\mathbb{R}^{n+1}\times\mathbb{R}_z$ by
\begin{equation}\label{eq:3.3}
\begin{array}{rl}
P_\psi v(z,t,x)=\int e^{i(x\cdot\xi+t\tau+z\sigma)}p(t,x,\tau,\xi+i|\sigma|\nabla\psi)\hat{v}(\sigma,\xi,\tau)d\sigma d\tau d\xi.
\end{array}
\end{equation}
It is very important to note that $p(t,x,\tau,\xi+i|\sigma|\nabla\psi)$ is independent of $t$ and $z$, which enable us to naturally define the pseudo-differential operator $P_\psi$.
In terms of the scaling
$$
(\xi,\sigma,\tau)\mapsto (\rho\xi,\rho\sigma,\rho^{2/\alpha}\tau)
$$
with large $\rho>0$.
We denote the principal symbol of $P_\psi$ by $\tilde{p}_\psi$ which is given by
\begin{equation}\label{SymbolTildePsi}
\begin{array}{rl}
\tilde{p}_\psi=&(1+i\tau)^{\alpha}+a_{nn}(\xi_n+i|\sigma|\tilde{X})^2\\
&+2\Sigma_{j=1}^{n-1}a_{jn}(\xi_j+2cx_j\xi_n+2icx_j|\sigma|\tilde{X})(\xi_n+i|\sigma|\tilde{X})\\
&+\Sigma_{j,k=1}^{n-1}a_{jk}(\xi_j+2cx_j\xi_n+2icx_j|\sigma|\tilde{X})(\xi_k+2cx_k\xi_n+2icx_k|\sigma|\tilde{X}),
\end{array}
\end{equation}
where $\tilde{X}=x_n-2X$ and $\alpha_1=\alpha$.
This $\tilde{p}_\psi$ coincides with the one given in \cite{LN2019} with the same notation.
Note that the form of the principal symbol $\tilde{p}_\psi$ is the same even for the case $2>\alpha=\alpha_1>\alpha_2>\cdots>\alpha_m>0$.

Recall the definition of Poisson bracket
\begin{equation}\label{3.4}
\{\Re \tilde{p}_\psi, \Im \tilde{p}_\psi\}=
\Sigma_{j=1}^{n}(\partial_{\xi_j}\Re \tilde{p}_\psi\cdot\partial_{x_j}\Im \tilde{p}_\psi-
\partial_{x_j}\Re \tilde{p}_\psi\cdot\partial_{\xi_j}\Im \tilde{p}_\psi)
+\partial_{\tau}\Re \tilde{p}_\psi\cdot\partial_{t}\Im \tilde{p}_\psi-
\partial_{t}\Re \tilde{p}_\psi\cdot\partial_{\tau}\Im \tilde{p}_\psi
\end{equation}
with the real part $\Re \tilde{p}_\psi$ and imaginary part $\Im \tilde{p}_\psi$ of $\tilde{p}_\psi.\,\,\text{Then}$ the principal part $\{\Re \tilde{p}_\psi, \Im \tilde{p}_\psi\}_p$ of the Poisson bracket
$\{\Re \tilde{p}_\psi, \Im \tilde{p}_\psi\}$ is
\begin{equation*}
\begin{array}{ll}
\{\Re \tilde{p}_\psi, \Im \tilde{p}_\psi\}_p=
\Sigma_{j=1}^{n}(\partial_{\xi_j}\Re \tilde{p}_\psi\cdot\partial_{x_j}\Im \tilde{p}_\psi
-\partial_{x_j}\Re \tilde{p}_\psi\cdot\partial_{\xi_j}\Im \tilde{p}_\psi).
\end{array}
\end{equation*}
Since $\tilde{p}_\psi$ are the same, we will have Lemma 3.2 and Lemma 4.1 of in  \cite{LN2019}. More precisely, we have the following estimate for the Poisson bracket and
sub-elliptic estimates for the operator $P_\psi$, respectively.
\begin{lemma}\label{lem3.2}
We have
\begin{equation}\label{3.8}
\begin{array}{l}
  (|\xi|^2+\sigma^2+|\tau|^{\alpha})^{3/2}\lesssim \{\Re \tilde{p}_\psi, \Im \tilde{p}_\psi\}_p
\end{array}
\end{equation}
when $\tilde{p}_{\psi}=0$ and \eqref{2.5} are satisfied.
\end{lemma}
\begin{lemma}\label{lem4.1}
There exists a small constant $z_0$
such that for all $u(t,x,z)\in C_0^{\infty}(U\times[-z_0,z_0])\cap\dot{\mathcal{S}}(\overline{{\mathbb{R}}_+^{1+n+1}})$, we have the estimates
\begin{equation}\label{4.1}
\begin{cases}
\begin{array}{l}
\Sigma_{k+s<2}||h(D_z)^{2-k-s}\Lambda_\alpha^{s}D_z^k u||+||h(D_z)\Lambda_\alpha u||\lesssim||P_\psi u||\quad {\rm if}\quad \alpha< \frac{4}{3},\\
\Sigma_{k+s<2}||h(D_z)^{2-k-s}\Lambda_\alpha^{s}D_z^k u||+||h(D_z)^{\frac{3}{2}-\frac{2}{\alpha}}\Lambda_\alpha^{\frac{2}{\alpha}} u||\lesssim||P_\psi u||\quad {\rm if}\quad \alpha\geq \frac{4}{3},
\end{array}
\end{cases}
\end{equation}
where $\dot{\mathcal{S}}(\overline{{\mathbb{R}}_+^{1+n+1}})$ is the subspace of Schwarz class $\mathcal{S}({\mathbb{R}}_+^{1+n+1})$ supported in $\{t\ge0\}$, $h(D_z)=(1+D_z^2)^{1/4}$, $U$ is an open neighborhood of the origin in ${\Bbb R}^{1+n}$ given as the interior of the set $\{(t,x): 0<t<T, |x'|\leq\mathcal{C}\sqrt{X}, 0\leq x_n\leq X\}\subset U$ with some constant $\mathcal{C}>0$,  and
$\Lambda_{\alpha}^s(\tau,\xi)=((1+|\xi|^2)^{1/\alpha}+i\tau)^{s\alpha/2}$.
\end{lemma}

Now by repeating the argument given for Lemma 5.1 of \cite{LN2019} which is based on Treve's argument, we can derive the following Carleman estimate for the operator $P(t,x,D_t,D_x)$.

\begin{lemma}\label{lem5.1} Let $\psi=\frac{1}{2}(x_n-2X)^2$.
Then there exists a sufficiently large constant $\beta_1>0$ depending on $n$
such that for any $\beta\geq \beta_1$, we have the estimates
\begin{equation}\label{5.1}
\begin{cases}
\begin{array}{l}
\sum_{|\gamma|\leq1}\beta^{3-2|\gamma|}\int e^{2\beta\psi(x)}|D_x^\gamma v|^2dtdx
\lesssim \int e^{2\beta\psi(x)}|P(t,x,D_t,D_x) v|^2dtdx,\quad \alpha<\frac{4}{3}\\
\sum_{|\gamma|\leq1}\beta^{3-2|\gamma|}\int e^{2\beta\psi(x)}|D_{x}^\gamma v|^2dtdx
+\beta^{3-\frac{4}{\alpha}}\int e^{2\beta\psi(x)}|D_{t}v|^2dtdx\\
\lesssim \int e^{2\beta\psi(x)}|P(t,x,D_t,D_x) v|^2dtdx,\quad \alpha\geq\frac{4}{3},
\end{array}
\end{cases}
\end{equation}
for any  $v(t,x)\in C_0^{\infty}(U)$.
\end{lemma}

Once having the Carleman estimates, we can show the conclusion of Proposition \ref{prop2.1} as follows.
First recall $u$ given just after \eqref{another 2.2} and the cutoff function $\chi$ given by \eqref{2.3}.
By an approximation argument using the Carleman estimate, we can assume that $u$ is smooth.
By \eqref{2.4},
$(\chi u)(t,x)\in C_0^{\infty}(U)$.
Then, by applying the Carleman estimates \eqref{5.1} to $(\chi u)(t,x)$, we have
\begin{equation}\label{6.2}
\begin{array}{rl}
&\sum_{|\gamma|\leq1}\beta^{3-2|\gamma|}\int_{x_n\leq (1-\varepsilon) X} e^{2\beta\psi(x)}|D^\gamma u|^2dtdx\\
\leq &\sum_{|\gamma|\leq1}\beta^{3-2|\gamma|}\int e^{2\beta\psi(x)}|D^\gamma (\chi u)|^2dtdx\\
\lesssim &\int e^{2\beta\psi(x)}|P(t,x,D_t,D_x) (\chi u)|^2dtdx\\
\lesssim &\sum_{|\gamma|\leq1}\int e^{2\beta\psi(x)}|\chi D^\gamma u|^2dtdx+
\int_{(1-\varepsilon) X<x_n\leq  X} e^{2\beta\psi(x)}|[P,\chi] u|^2dtdx,
\end{array}
\end{equation}
where $[\cdot,\cdot]$ denotes the commutator. Let $\beta$ be large enough to absorb the first term on the right hand side of \eqref{6.2} into the first and second lines. Then we have
\begin{equation}\label{6.3}
\begin{array}{rl}
&\beta^{3}\int_{x_n\leq (1-2\varepsilon) X} e^{\beta (1+2\varepsilon)^2 X^2/2}|u|^2dtdx\\
\leq &\sum_{|\alpha|\leq1}\int_{(1-\varepsilon) X<x_n\leq  X} e^{2\beta\psi(x)}|D^\alpha u|^2dtdx\\
\lesssim &C(u)e^{\beta (1+\varepsilon)^2 X^2/2}.
\end{array}
\end{equation}
Now let $\beta$ tend to $\infty$, then we obtain that $u=0$ on $x_n\leq (1-2\varepsilon) X$.
Since this argument works for any $1>\varepsilon>0$, we get that
$u=0$ in $\{x_n< X\}$ which implies the conclusion of Proposition \ref{prop2.1}.
\eproof

\section{Proof of Theorem \ref{thm1.3}}\label{sec7}

Let $z \in \partial B_r$. Then, from Proposition \ref{prop2.1}, there exists a small $0<\tilde{r}<1$ such that
$u=0$ in $(T_1,T_2)\times B_{\tilde{r}}(z)$ with $0<T_1<T_2<T$. Thus, we have the following equations.
\begin{equation*}
\left\{
\begin{array}{l}
\Sigma_{j=1}^mq_j\partial_t^{\alpha_j} u(t,y)-\emph{L} u(t,y)=\tilde{\ell}_1(t,y;\nabla_{y}u(t,y)),\\
u(t,y)=0,\,t\in(T_1,T_2),\,y\in B_{\tilde{r}}(z),\\
u(t,y)=0,\,t\in(0,T),\,y\in B_{r}.
\end{array}
\right.
\end{equation*}
Define $$R_r=\{(y_1,\cdots,y_n): |y_j|<r, \,  1\leq j\leq n\}$$ and $$R_r^-=\{y=(y_1,\cdots,y_n):y\in R_r,\, y_n<0, \}.$$
If necessary, rotate
the coordinates such that the line passing through the origin and $z\in\partial B_r$ is on the $y_n$ axis. Note that this will not affect \eqref{1.2}. By a diffeomorphism from $B_r\cup B_{\tilde r}(z)$ onto a neighborhood of  $\overline{R_1\cup(\hat{z}+R_{r_1})}$ with $\hat z=(0,\cdots,0,1)$ and some $0<r_1<1$, the push forward of $u$ satisfies $u=0$ in  $(0,T)\times R_1$ and $u=0$ in  $(T_1,T_2)\times \{\hat{z}:\hat{z}=(0,\cdots,1)+w,w\in R_{r_1}\}$, where we have abused the notations $y$ and $u$ to denote the push forward $u(t,y)$ of $u$. This diffeomorphism may change the constant $\delta$ in \eqref{1.2}, but we will use the same $\delta$ for the new \eqref{1.2}.
We also move $(0,\cdots,1)$ to the origin and use the same coordinates $y$ so that we have
$u=0$ in  $(0,T)\times R_{r_1}^-$ and $u=0$ in  $(T_1,T_2)\times R_{r_1}$. Further, we transform $u$ under the spatial diffeomorphism from $R_{r_1}$ to $R_1$ so that we have $u=0$ in  $(0,T)\times R_1^-$ and $u=0$ in  $(T_1,T_2)\times R_1$. Again, we will use the same notations $y$ for $u(t,y)$ and $\delta$ for \eqref{1.2}.
Thus $u$ satisfies the following equations.
\begin{equation*}
\left\{
\begin{array}{l}
\Sigma_{j=1}^mq_j\partial_t^{\alpha_j} u(t,y)-\emph{L} u(t,y)=\tilde{\ell}_1(t,y;\nabla_{y}u(t,y)),\\
u(t,y)=0,\,t\in(T_1,T_2),\,y\in R_1,\\
u(t,y)=0,\,t\in(0,T),\,y\in R_1^-.
\end{array}
\right.
\end{equation*}

Now we consider the diffeomorphism $\tilde{y}$ from $R_1$ to ${\mathbb R}^n$ by
\begin{equation}\label{diffeomorphism}
\tilde{y}(y)=(\frac{y_1}{\sqrt{1-y_1^2}},\cdots,\frac{y_n}{\sqrt{1-y_n^2}}).
\end{equation}
It is not hard to see that
$$\tilde{\emph{L}} u(t,\tilde{y}):=\emph{L} u(t,y)=\sum_{j,k=1}^n\tilde{a}_{j,k}(t,\tilde{y})\partial_{\tilde{y}_j}\partial_{\tilde{y}_k}u(t,\tilde{y})+{l}_1(t,\tilde{y};\nabla_y)u(t,\tilde{y}),$$
where $\tilde{a}_{j,k}(t,\tilde{y})=a_{j,k}(t, {y})(1+\tilde{y}_j^2)^{3/2}(1+\tilde{y}_k^2)^{3/2}$. From \eqref{1.2}, we have $\tilde{a}_{j,k}(t,\tilde{y})$ satisfying
\begin{equation}\label{3.2}
\begin{array}{l}
\delta |\tilde{\eta}|^2 \leq \sum_{j,k=1}^n\tilde{a}_{j,k}(t,\tilde{y})\eta_j\eta_k=\sum_{j,k=1}^n a_{j,k}(t,\tilde{y})\tilde\eta_j\tilde\eta_k \leq \delta^{-1} |\tilde{\eta}|^2,
\end{array}
\end{equation}
where $\tilde{\eta}=((1+\tilde{y}_1^2)^{3/2}\eta_1,\cdots,(1+\tilde{y}_n^2)^{3/2}\eta_n)$
and $\eta_j$ corresponds to the symbol of $D_j:=-i\partial_{\tilde{y}_j}$.
Thus, by this diffeomorphism, $R_1$ is mapped to  ${\mathbb R}^n$ and $u(t,\tilde y)$ satisfies the following equations.
\begin{equation}\label{3.3}
\left\{
\begin{array}{l}
\Sigma_{j=1}^mq_j\partial_t^{\alpha_j} u(t,\tilde y)-\tilde{\emph{L}} u(t,\tilde{y})=\tilde{\ell}_1(t,\tilde{y};\nabla_{\tilde y})u(t,\tilde{y}),\\
u(t,\tilde{y})=0,\,t\in(T_1,T_2),\,\tilde{y}\in {\mathbb R}^n,\\
u(t,\tilde{y})=0,\,t\in(0,T),\,\tilde{y}\in {\mathbb R}_-^n.
\end{array}
\right.
\end{equation}

\medskip
Then we first have the following proposition.
\begin{pr}\label{Pro3.1}
Assume that  $u(t,\tilde{y})\in H^{\alpha,2}((0,T)\times\mathbb{R}^n)$ satisfies \eqref{3.3}.
Then $u(t,\tilde{y})=0$ in $(T_1,T)\times\mathbb{R}^n$.
\end{pr}
\pf.

Let $\theta$ be a smooth function such that
\begin{equation}\label{3.4}
\theta (t)=
\begin{cases}
\begin{array}{l}
0,\quad t\leq T_1,\\
1,\quad t\geq T_2.
\end{array}
\end{cases}
\end{equation}
Let $\hat{u}:=\theta (t)u$ and compute $\partial_{t}^\alpha \hat{u}$ as follows.
\begin{equation}\label{3.5}
\begin{array}{l}
\Gamma(1-\alpha)\partial_{t}^\alpha \hat{u}=\int_0^t(t-\eta)^{-\alpha}\partial_\eta (\theta(\eta)u(\eta,y))d\eta\\
=
\begin{cases}
\begin{array}{l}
0,\quad t\leq T_1,\\
0,\quad T_1<t< T_2.\\
\int_{T_2}^{t}(t-\eta)^{-\alpha}\theta(\eta)\partial_\eta u(\eta,y)d\eta=\int_0^t(t-\eta)^{-\alpha}\theta(\eta)\partial_\eta u(\eta,y)d\eta,\quad T_2\leq t.
\end{array}
\end{cases}\\
=\theta(t)\int_0^t(t-\eta)^{-\alpha}\partial_\eta u(\eta,y)d\eta\\
=\Gamma(1-\alpha)\theta(t)\partial_{t}^\alpha u.
\end{array}
\end{equation}
Then, by multiplying $\theta(t)$ to the each equation in \eqref{3.3}, we have
\begin{equation}\label{3.6}
\left\{
\begin{array}{l}
\Sigma_{j=1}^m q_j\,\partial_t^{\alpha_j} \hat{u}(t,\tilde y)-\tilde{\emph{L}} \hat{u}(t,\tilde{y})=\tilde{\ell}_1(t,\tilde{y};\nabla_{\tilde y})\hat{u}(t,\tilde{y}),\\
\hat{u}(t,\tilde{y})=0,\,t\leq 0,\,\tilde{y}\in \mathbb{R}^n,\\
\hat{u}(t,\tilde{y})=0,\,t\in(0,T),\,\tilde{y}\in \mathbb{R}_-^n.
\end{array}
\right.
\end{equation}
Now, we apply Theorem 7.3 of \cite{LN2019} to get $\hat{u}=\theta (t)u=0$ for $t>0$.
By the definition of $\theta$, we have that $u(t,\tilde{y})=0$ for $t>T_1$.
\eproof

In order to complete the proof of Theorem \ref{thm1.3}, we will prove the
following proposition.
\begin{pr}\label{Pro3.2}
Assume that  $u(t,\tilde{y})\in H^{\alpha,2}((0,T)\times\mathbb{R}^n)$ satisfies \eqref{3.3}.
Then $u(t,\tilde{y})=0$ in $(0,T)\times (B'_{r_2}\times\mathbb{R})$, where $B'_{r_2}\subset\mathbb{R}^{n-1}$ is a ball centered at the origin with a small radius $r_2>0$.
\end{pr}
\pf.
From Proposition \ref{Pro3.1}, we can assume that $u(t,\tilde{y})=0$ for $t>T_1$.
For the rest, it suffices to prove that
\begin{equation}\label{3.7}
\begin{array}{l}
u(t,\tilde{y})=0\quad{\rm in} \quad \displaystyle (0,T_1)\times(B'_X\times{\mathbb R})=\cup_{s\in \mathbb{N}} E_s,
\end{array}
\end{equation}
where
$$E_s=\{(t,\tilde{y}):t\in(0,T_1],\, \tilde{y}_n+|\tilde{y}'|^2-\frac{sX}{T_1}t<0,\, |\tilde{y}'|<X\}$$
with a fixed small positive constant $X$. The idea of the proof is the same as that of Theorem 7.3 of \cite{LN2019}.

Let
\begin{equation*}
\begin{array}{l}
\begin{cases}
Q_s=\{(\tilde{y}_1,\cdots,\tilde{y}_n):|\tilde{y}'|\leq X,-1/6<\tilde{y}_n\leq  sX \},\\
\tilde{Q}_s=\{(\tilde{y}_1,\cdots,\tilde{y}_n):|\tilde{y}'|\leq 2X,-1/3<\tilde{y}_n\leq  (1+s)X\},
\end{cases}
\end{array}
\end{equation*}
where $\tilde{y}'=(\tilde{y}_1,\cdots,\tilde{y}_{n-1})$. Then define a cut off function
$\kappa_s(\tilde{y})\in C^\infty(\mathbb{R}^n)$ such that \begin{equation*}
\begin{array}{l}
\kappa_s(\tilde{y})=
\begin{cases}
1, \quad \tilde{y}\in Q_s,\\
0, \quad \tilde{y}\in \mathbb{R}^n\setminus \tilde{Q}_s.
\end{cases}
\end{array}
\end{equation*}

\medskip
We first show that $u=0$ in $E_1$. For this, Figure 1 given below will be helpful for understanding the coming argument. Note that we clearly have
$$(\kappa_1 u)(t,\tilde{y})=0,\quad \partial^{\alpha}_{t}(\kappa_1 u)(t,\tilde{y})=0 \quad {\rm for} \quad \tilde{y}_n \leq 0\quad  {\rm or}\quad \tilde{y}\in \mathbb{R}^n\setminus \tilde{Q}_1$$
and
$$(\kappa_1 u)(t,\tilde{y})=u(t,\tilde{y}), \quad \tilde{y}\in Q_1.$$

Even though the proof is the same as that of Theorem 7.3 of \cite{LN2019}, we need to modify the Holmgren type transform as follows.
\begin{equation}\label{3.8}
x'=\tilde{y}', x_{n}=\tilde{y}_n+|\tilde{y}'|^2-\frac{X}{T_1}t+X, \tilde{t}=t.
\end{equation}
We call this coordinates the {\sl Holmgren type coordinates}. It is clear that
\begin{equation}\label{3.9}
{\rm supp}\, \chi\kappa_1 u\subset \tilde{Q}_1\cap\{0\leq x_n< X\},
\end{equation}
where $\chi$ is the one define by \eqref{2.3}.

 \begin{figure}[h!]\label{fig1}
	\centering
	\includegraphics[scale=0.35]{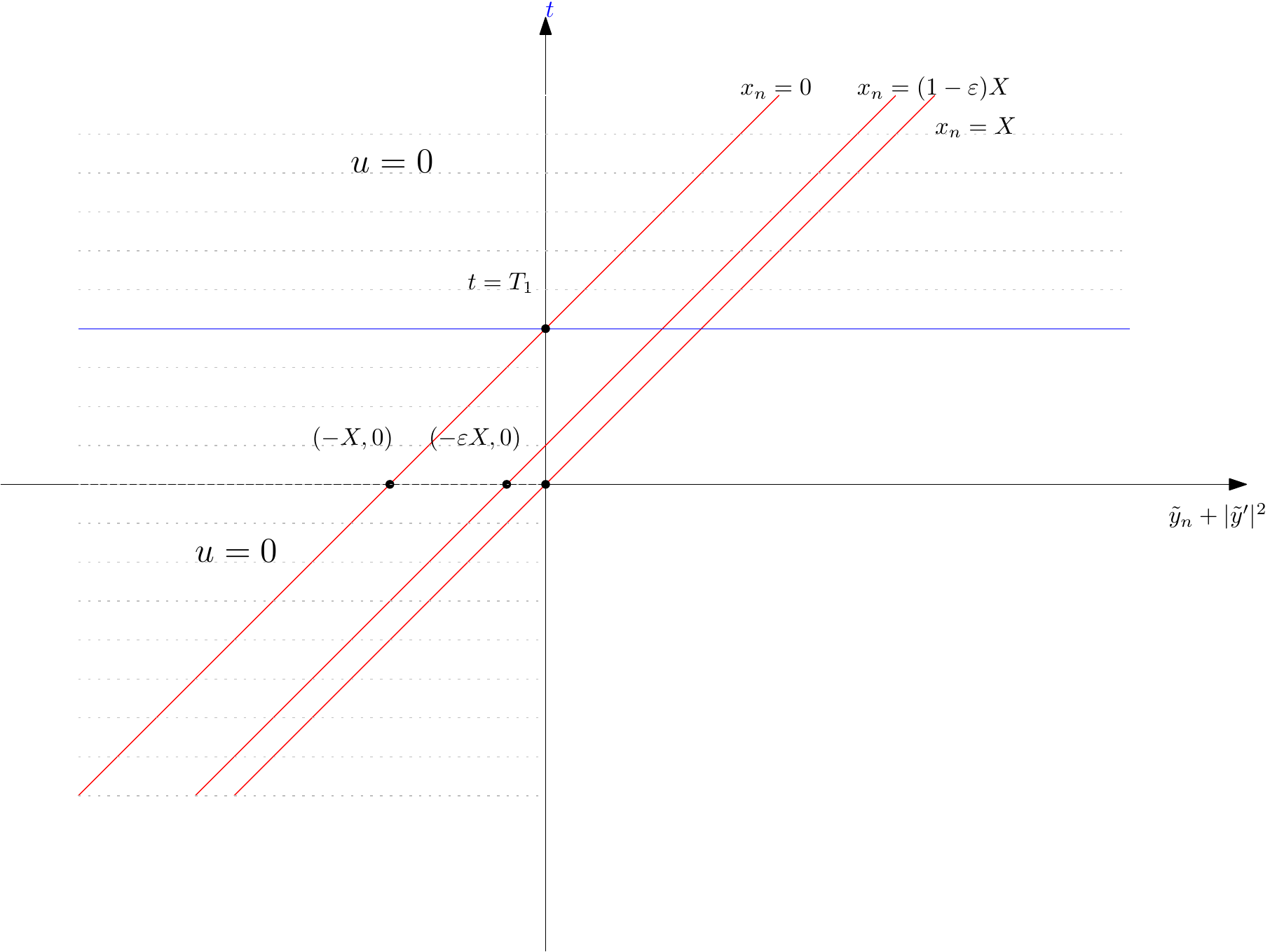}
	\caption{$x_{n}=\tilde{y}_n+|\tilde{y}'|^2-\frac{X}{T_1}t+X$.}
\end{figure}

Since $\kappa_1 u=u$ in $Q_1$, we will only denote $\kappa_1 u$ by $u$. Also, we will assume that
\begin{equation*}
|x'|\lesssim X<\sqrt{X}, \quad x_n\leq X.
\end{equation*}
As before, upon conjugating $\sum_{l=1}^m q_l\,\partial_t^{\alpha_l}-\tilde{\emph{L}}$ with
$e^{t}$, we denote what we get by $P(t,x,D_t,D_x)$
in terms of the Holmgren type coordinates $(t,x)$ (see \eqref{3.8}). Then it is easy to see that the principal symbol $\tilde p_\psi$ of $P_\psi=P_\psi(t,x,D_t,D_{x},D_z):=P(t,x,D_t,D_x+i|D_z|\nabla\psi)$ with the Carleman weight function $\psi=\frac{1}{2}(x_n-2X)^2$ is given as
\begin{equation}\label{3.10}
\begin{array}{rl}
\tilde p_\psi=&\Sigma_{l=1}^m q_l(1+i\tau)^{\alpha_l}+\tilde a_{nn}(t,x)(\xi_n+i|\sigma|\tilde X)^2\\
&+2\Sigma_{j=1}^{n-1}\tilde a_{jn}(t,x)(\xi_n+i|\sigma|\tilde X)(\xi_j+2cx_j\xi_n+2icx_j|\sigma|\tilde X)\\
&+\Sigma_{j,k=1}^{n-1}\tilde a_{jk}(t,x)(\xi_j+2cx_j\xi_n2ic x_j|\sigma|\tilde X)(\xi_k+2cx_k\xi_n+2icx_k|\sigma|\tilde X)
\end{array}
\end{equation}
with $\tilde X=x_n-2X$. Hence this equation for $\tilde p_\psi$ coincide with \eqref{SymbolTildePsi} except the difference in notations. Then, based on this, we can argue as in the proof of Theorem 7.3 in \cite{LN2019} to show $u(t,\tilde{y})=0$ in $E_1$.

Next we show that $u=0$ in $E_s$ for $s\ge 2$. For this, Figure 2 given below will be halpful for understanding the situation we have in the coming argument. When $s=2$, we need to take the Holmgren type transform as follows.
\begin{equation}\label{3.11}
x'=\tilde{y}', x_{n}=\tilde{y}_n+|\tilde{y}'|^2-\frac{2X}{T_1}t+X, \tilde{t}=t.
\end{equation}
Then by observing
\begin{equation}\label{3.12}
{\rm supp}\, \chi\kappa_2 u\subset \tilde{Q}_2\cap\{0\leq x_n< X\},
\end{equation}
we have $u(t,\tilde{y})=0$ in $E_2$ with
$$E_2=\{(t,\tilde{y}):t\in(0,T_1],\, \tilde{y}_n+|\tilde{y}'|^2-\frac{2X}{T_1}t<0,\, |\tilde{y}'|<X\}$$
by repeating the same argument as we did for $s=1$.
Now, we assume that $u(t,\tilde{y})=0$ in $E_s$ with
$$E_s=\{(t,\tilde{y}):t\in(0,T_1],\, \tilde{y}_n+|\tilde{y}'|^2-\frac{sX}{T_1}t<0,\, |\tilde{y}'|<X\}.$$
Then, using the Holmgren type transform of the form
\begin{equation}\label{3.13}
x'=\tilde{y}', x_{n}=\tilde{y}_n+|\tilde{y}'|^2-\frac{(s+1)X}{T_1}t+X, \tilde{t}=t
\end{equation}
and observing
\begin{equation}\label{3.14}
{\rm supp}\, \chi\kappa_{s+1} u\subset \tilde{Q}_{s+1}\cap\{0\leq x_n< X\},
\end{equation}
we have $u(t,\tilde{y})=0$ in $E_{s+1}$ with
$$E_{s+1}=\{(t,\tilde{y}):t\in(0,T_1],\, \tilde{y}_n+|\tilde{y}'|^2-\frac{(s+1)X}{T_1}t<0,\, |\tilde{y}'|<X\}$$
by repeating the same argument. In this way we can obtain \eqref{3.7}. Then, combining this with Proposition \ref{Pro3.1}, we immediately have the conclusion of Proposition \ref{Pro3.2}.

 \begin{figure}[h!]\label{fig2}
	\centering
	\includegraphics[scale=0.35]{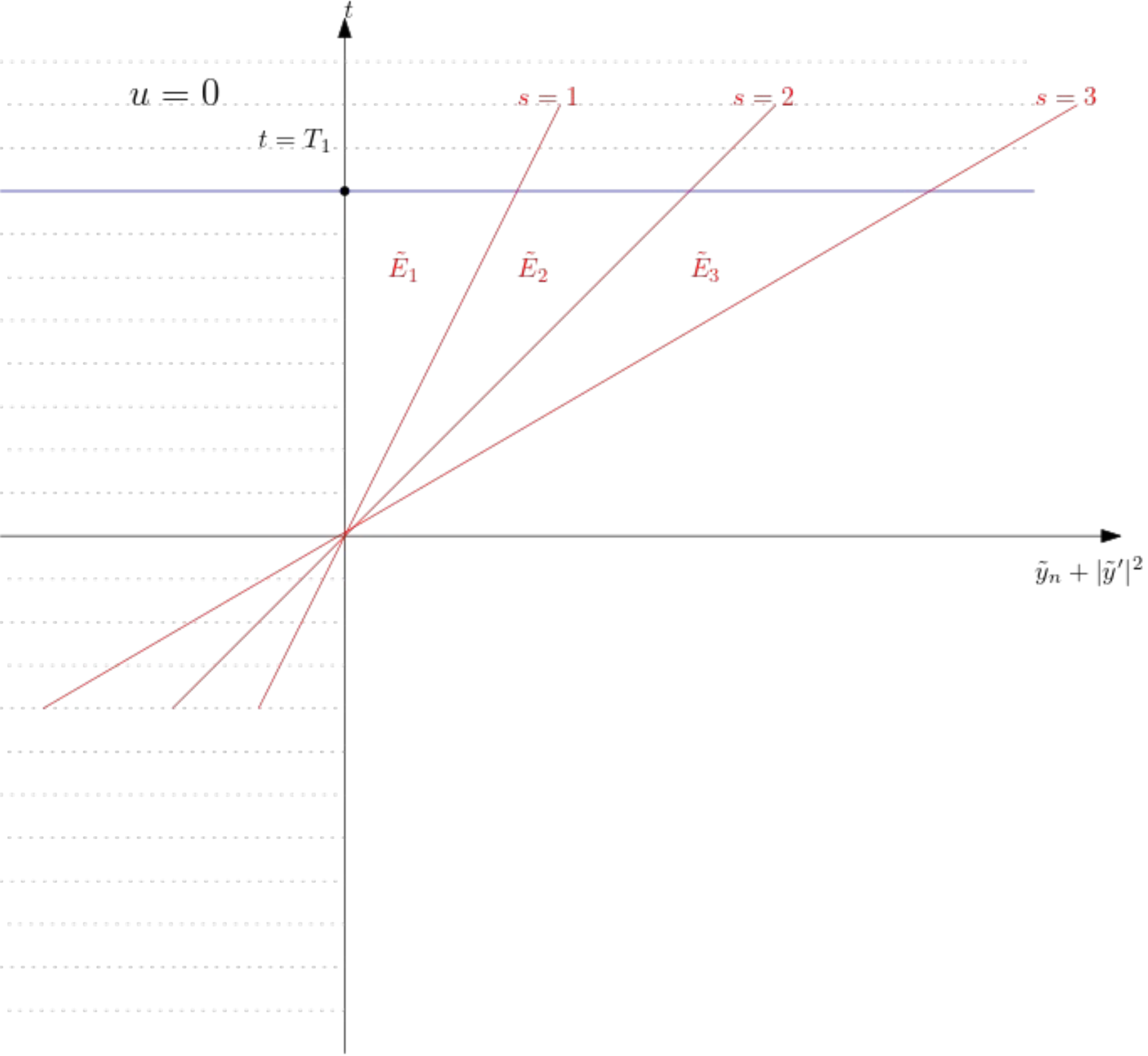}
	\caption{The red lines are $\tilde{y}_n+|\tilde{y}'|^2=\frac{sX}{T_1}t$ and  $E_{s}=\cup_{j=1}^s\tilde{E}_j$.}
\end{figure}

\eproof

\bigskip
\noindent
{\bf \large  Acknowlegement}${}$
\newline
The first author was partially supported by the
Ministry of Science and Technology of Taiwan. Also, the second author acknowledge the support by Grant-in-Aid for Scientific Research of the Japan Society for the Promotion of Science (No.19K03554) during this study.

\end{document}